%% file: StabCont.tex
\documentclass[a4paper, 11 pt, peerreview]{ieeeconf}  

\IEEEoverridecommandlockouts                              

\usepackage{graphics} 
\usepackage{epsfig} 
\usepackage{amsmath} 
\usepackage{amssymb}  

\usepackage{amsfonts}
\usepackage{color}
\usepackage{txfonts}
\usepackage{pifont}
\usepackage{url}
\usepackage{flushend}
\usepackage{hyperref}
\usepackage[normalem]{ulem}

\newcommand{\calW}{\mathcal{W}}
\newcommand{\RR}{\mathbb{R}}

\newcommand{\ud}{\mathrm{d}}
\newcommand{\T}{\mathrm{T}}

\title{\LARGE \bf
Continuation for stability domain determination: a case study
}

\author{Quentin Peyron$^{1}$, Isabelle Charpentier$^{2}$ and Edouard Laroche$^{2}$
\thanks{$^{1}$Quentin Peyron is with the Institut National des Sciences Appliqu\'ees, F-67084 Strasbourg Cedex, France 
        {\tt\small quentin.peyron@insa-strasbourg.fr}}%
\thanks{$^{2}$Isabelle Charpentier and Edouard Laroche are with the ICube laboratory, University of Strasbourg and CNRS, 67412 Illkirch Cedex, France 
        {\tt\small \{icharpentier,laroche\}@unistra.fr}}%
}

\begin{document}

\maketitle
\thispagestyle{empty}
\pagestyle{empty}

\begin{abstract}
Determining a stability domain, i.e. a set of equilibria for which a dynamical system remains stable, is a core problem in control. When dealing with controlled systems, the problem is generally transformed into a robustness analysis problem: considering a set of uncertain parameters (generally chosen as polytopic), one seeks for guaranties that the system remains stable for this given set or some dilatation or contraction of this set. Thanks to the development of formal methods based on continuation, another approach is proposed in this paper to determine the border of the stability domain. The proposed developments were made with MatCont, a package to be used with Matlab, and the considered system is a 3-cable robot with linear PD and PID controllers. 
\end{abstract}

\section{INTRODUCTION}

Determining the stability of an equilibrium point of a dynamic system is a very central problem in control. For systems modeled by ordinary differential equations (ODE), it is well known that the local stability is equivalent to the Jacobian matrix of the state function being Hurwitz, i.e. it must have all its eigenvalues located in the left-hand half plane. 

Let us consider the model 
\begin{equation} \label{eq:ode}
 \dot{\chi} = f(\chi,\alpha),
\end{equation}
where the state vector $\chi\in\RR^n$ and the vector-valued function $f:\RR^n\times \RR^p \rightarrow \RR^n$ depend on the variations of the parameter vector $\alpha\in\RR^p$. Determining the set of parameters and equilibrium states for which the system remains stable is a much more complex issue than simply proving the stability for one particular value of $\alpha$. For (\ref{eq:ode}), an equilibrium $\chi_e$ verifies $f(\chi_e,\alpha) = 0$ and is locally stable if the Jacobian matrix $f_\chi(\chi_e,\alpha) = \frac{\ud f}{\ud \chi}(\chi_e,\alpha)$ has no eigenvalue with positive real part. 

Notice that in control, many contributions have been devoted to solving the robustness analysis problem \cite{Packard93,FAG96}. Considering a nominal value $\alpha_0$ for which the system is stable and a set $E_\alpha$ including $\alpha_0$, the aim is to prove that the system is stable for any $\alpha \in E_\alpha$. Considering the set $\rho E_\alpha$ obtained from a dilatation of $E_\alpha$ of ratio $\rho$, the largest value $\rho^*$ such that the system is stable in $\rho^* E_\alpha$, is the robustness margin. Determining accurately the robustness margin in an efficient manner is a difficult problem. 

The objective of the paper is to test the potentiality of continuation methods for determining the limit of the stability domain. Starting from a stable equilibrium $\chi_e(\alpha_0)$, variations of $\alpha$ are first performed according to one of its components until an eigenvalue has its real part equals zero, i.e. reaching a limit of the stability domain. In a second time, two directions can be chosen in the parameter space to track the equilibrium, yielding a closed 2D curve. Notice that this approach does not guaranty stability for all points inside the obtained curve. Additional explorations in a number of directions can detect possible inner unstable regions. The purpose of this paper is to present the use of a numerical package for continuation analysis, MatCont, for the evaluation of the stability domain. To our best knowledge, the use of continuation for the determination of stability domains has not been considered yet. 

Cable-driven parallel robots (CDPR) have been the focus of research interests for the last decade as can attest the two conferences dedicated on the field \cite{BP13,BP15}. Composed of a platform moved by cables that are rolled at the distant end, these robots may be adequate solutions for a number of issues thanks to their low mass, large workspace and low invasiveness. However, a number of issues are to be solved in order to obtain accurate positioning. The evaluation of the continuation for stability analysis for CDPR has been done on a simplified case which includes nonlinearities but remains of reasonable complexity in order to yield tractable computations. 

The paper is organized as follows. Section~2 introduces the continuation methods and the MatCont package. Section~3 details the cable robot, provides its model and the considered control laws. Section~4 reports the numerical results, i.e. the stability domains obtained with MatCont. A conclusion is provided in Section~5. 

\section{CONTINUATION METHODS FOR STABILITY DOMAIN DETERMINATION}\label{sec:continuation}

On the one hand determining the stability of an equilibrium of a dynamic system is a very central problem in control. On the other hand, continuation has become a classical tool for bifurcation and stability analysis of ODE (\ref{eq:ode}), differential-algebraic equations (\ref{eq:algebraic}) or partial differential equations,
\begin{equation} 
f(\chi,\alpha)=0, \label{eq:algebraic}
\end{equation}
in research domains ranging from  mechanics to economics \cite{charpentier12}. This section introduces basics of continuation  together with a brief state-of-the-art in robotics considering (\ref{eq:ode}) and (\ref{eq:algebraic}), respectively. MatCont's abilities for the determination of stability domains are then presented.
 
\subsection{Basics of continuation and bifurcation analysis}

Let $\lambda$ be a scalar component of the parameter set $\alpha=\hat{\alpha}\cup\{\lambda\}$ and $\hat{\alpha} = \alpha\setminus\{\lambda\}$ be the other ones. Continuation is a solution method for (\ref{eq:ode}), for instance, with $\lambda$ to be varied. The solutions $X = (\chi,\lambda)$
thus form branches of solutions. Bifurcation and stability analysis of nonlinear problems is one of the major pillars in computational sciences, see \cite{Kel87,AG90,MatCont,Sey09} and the references therein. The pseudo-arc length continuation method \cite{Kel77} is of particular interest since it introduces a path parameter (measuring the pseudo-arc length along the branch) to close the under-determined nonlinear problem under study. From a numerical point of view, continuation is carried out using either a first order Newton-Raphson method \cite{MatCont,DKK91} or a higher order one \cite{CDPF07,charpentier12}. At a given solution point $X$, these iterative predictor-corrector methods require the calculation of the Jacobian $f_X(X)$ and the tangent vector $V(X)$ to follow the branch of solutions. 

These two derivatives also serve in the detection of bifurcations (singular points) such as folds (also called limit points, LP), Hopf bifurcations (H), and branch points (BP) where two solution branches intersect. Their common characteristics is a null singular value for $f_X$. Test functions allow to classify them. For instance, the BP test monitors any change of sign of the determinant 
\begin{equation}
\Delta(X)=\det\left(\begin{bmatrix} f_X(X) & V(X) \end{bmatrix} \right).
\end{equation} A fold is characterized by $\Delta(X)\ne 0$ and a null value for the component of $V(X)$ related to $\lambda$. A Hopf bifurcation corresponds to a steady point of (\ref{eq:ode}) for which a local change in the stability properties results in the appearance or the disappearance of a limit cycle in the time-domain trajectories. As described in Section \ref{sec:results}, bifurcation points are of particular interest in the study of workspace boundaries and stability domains.

\subsection{Determination of workspace boundaries}\label{sec:workspaceboundaries}

Mechanism equations may be frequently written in the general form (\ref{eq:algebraic}), where the vectors of input and output coordinates are denoted by $\chi$ and $\alpha$, respectively. The reachable workspace $\mathcal{W}$ may be defined as
\begin{equation}
\mathcal{W} = \{ \alpha | ~ \exists ~ \chi ~ \text{such that} ~ f (\chi,\alpha) = 0 \}.
\end{equation}
The workspace boundaries $\partial \mathcal{W}$, see \cite{litvin80,haug96,crmeca}, 
\begin{equation}
\partial \mathcal{W} = \{ \alpha \in \mathcal{W} ~ | ~ \exists ~ (\chi,\zeta) ~ \text{s.t.} ~ f_{\chi}^\T(\chi,\alpha) \cdot \zeta = 0 \hbox{~and~} \zeta^\T \cdot \zeta=1 \},
\label{eq:boundaryCondition}
\end{equation}
correspond to the singular configurations of the mechanism, with $\zeta$ a normal vector to $\partial \mathcal{W}$ and $\zeta^\T$ denoting the transpose of $\zeta$. They may be determined by applying a continuation method to the extended system of equations \cite{haug96}
\begin{equation}\label{eq:extended}
f^{\partial \calW}(\chi,\alpha,\zeta)=
\left(\begin{array}{c}
	f(\chi,\alpha)\\
	f_\chi^\T(\chi,\alpha) \cdot \zeta\\
	\zeta^\T \cdot \zeta -1
\end{array}\right)
= 0.
\end{equation}
It is worth noticing that this system is very similar to the extended systems described in \cite{Sey79,Sey09} and used in MatCont \cite{MatCont}. 

\subsection{Dynamical system and equilibrium curve in MatCont}\label{sec:MatCont}

MatCont is a Matlab software project developed by Govaerts and Kuznetsov, and devoted to the numerical continuation and bifurcation studies of parameterized dynamical systems (\ref{eq:ode}). MatCont is freely \href{http://www.MatCont.ugent.be/}{available} for non-commercial use. The software abilities are described with details in the MatCont manual \cite{MatCont}. Fig. \ref{fig:MatCont} illustrates some of the MatCont's functionalities we use to implement the stability domain determination of the CDPR models presented in Section \ref{sec:CDPR}.

\begin{figure}
\begin{center}
\includegraphics[width=3.8cm]{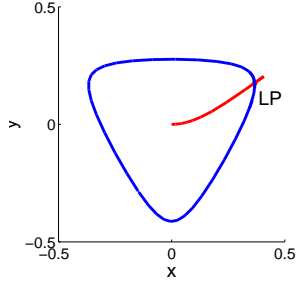}
\caption{Continuation process. Red: Solution for (\ref{eq:ode}). Blue:  Fold curve and stability region obtained by solving  (\ref{eq:algebraic}) from the limit point (LP)}\label{fig:MatCont}
\end{center}
\end{figure}
MatCont allows for the solution of the dynamical system (\ref{eq:ode}) with respect to both the time and some modeling parameter $\lambda$, performing a so-called one-codim continuation. One of the solution branches $X=(\chi,\lambda)$ is plotted as a red line in the projected bifurcation diagram Fig \ref{fig:MatCont}. Bifurcation tests are performed along the continuation process in the search for singular points on the computed branch. In Fig.~\ref{fig:MatCont}, a fold (LP) is detected at some point $X_e = (\chi_e,\alpha_e)$. This singular point satisfies (\ref{eq:algebraic}) as well. An equivalent criterion \cite{Sey79,Sey09} is that there exist a nontrivial $\zeta$ such that the linearized equation,
\begin{equation}
f_{(\chi,\hat{\alpha})}(\chi,\hat{\alpha},\lambda) \cdot \zeta=0,
\end{equation}
is satisfied. This allows to build an extended system very similar to (\ref{eq:extended}). If the initial point $(\chi_0,\alpha_0)$ chosen to solve \eqref{eq:ode} is a stable point, then the branch (excepted LP) running from  $(\chi_0,\alpha_0)$ (at point (0,0) in Fig. \ref{fig:MatCont}) to the LP point is a stable one. It is worth noting that MatCont allows for the calculation of unstable branches, i.e. after the LP point, what a standard method based on ODE cannot do. At this LP point, MatCont allows to build a fold curve (blue line in Fig. \ref{fig:MatCont}) through a two-codim continuation, i.e. with respect to two parameters. 

From a computer point of view, MatCont operates the continuation from the user-defined equations. Jacobian calculation as well as the extended system generation are hidden to the user \cite{MatCont}. Continuation may be carried out using either the MatCont's Graphical User Interface or the CL\_MatCont command lines. One-dimensional solution branches are computed very accurately. 

For the sake of completeness, two general limitations are reported. First, continuation computes 1D curves. For higher-dimensional studies, projections in several planes must be considered. Second, the solution branch issued from $(\chi_0,\alpha_0)$ can only detect instability regions it intersects. Options to get a global result range from a first global analysis ($\mu$-analysis) or the choice of different initial points or different continuation parameters to run complementary simulations. Notice that these limitations are inherent to any continuation approach and are not specific to MatCont. 

\section{STUDY CASE: A CABLE-DRIVEN PARALLEL ROBOT }\label{sec:CDPR}

\begin{figure}
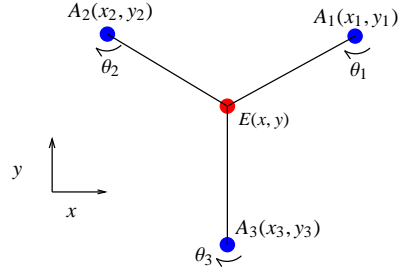

\begin{center}
\input robot_pspdf.tex
\caption{Geometry of the CDPR} 
\label{fig:schemaR3CR}
\end{center}
\end{figure}

The system, depicted in Fig.~\ref{fig:schemaR3CR}, has been chosen as simple enough for evaluating the efficiency of the continuation approach. 

\subsection{System description and model }\label{model}

The platform is assimilated to a punctual mass $m$ moving in the (O, $x$, $y$) horizontal plane thanks to three cables rolled at fixed ends A$_k$, $k=1,2,3$. The coordinates of the platform are ($x(t)$, $y(t)$). The central location corresponds to $x=0$ and $y=0$. Roller \#$k$ of radius $R$, of angular position denoted $\theta_k(t)$, with $\theta_k = 0$ corresponding to a fully rolled cable, is actuated by a torque $\Gamma_k$ with infinite bandwidth. The inertias of the rollers are assumed constant equal to $I$, thus neglecting the contribution of the winded cables. Cables are assumed to be strait without elasticity nor mass. A number of models, possibly with more cables and flexible cables are available in the literature, see for instance \cite{KZW06,LCCG12}. 

Denoting $q=\begin{bmatrix} x & y \end{bmatrix}^\T$, the inverse kinematic model writes
\begin{equation}
  \theta_k(q) = \frac{1}{R} \sqrt{(x - x_k)^2 + (y - y_k)^2}, \quad\hbox{for }k=1,..,3,
\label{MGI}
\end{equation}
and the differential kinematic model is given by
\begin{equation}
  \dot\theta_k(q,\dot q) = J(q) \, \dot q,
\label{MCI}
\end{equation}
where $J(q) = \begin{bmatrix} \frac{\partial \theta}{\partial x}(q) & \frac{\partial \theta}{\partial y}(q) \end{bmatrix}$.

The dynamic model can be written either with Euler-Lagrange or Newton-Euler approaches. The first one is considered herein. The Lagrangian includes only kinetic energy and writes
\begin{equation}
  L(q,\dot{q}) = \frac{1}{2} m \, (\dot{x}^2+\dot{y}^2) + \frac{1}{2} I \, (\dot{\theta_1}^2 +\dot{\theta_2}^2 + \dot{\theta_3}^2).
\end{equation}
This kinetic energy can be re-witten $L(q,\dot{q}) = \frac{1}{2} \dot q^\T \, M(q) \, \dot q$ where the inertia matrix $M(q)$ can be easily computed. 

Euler-Lagrange equations then provide the dynamic model
\begin{equation}
	M(q) \, \ddot{q} + C(q,\dot{q}) \, \dot{q} = J^{\T}(q,\dot{q}) \, \Gamma ,
	\label{équations non linéaires}
\end{equation}
where $C(q,\dot{q})$ includes the Coriolis effects plus $b$, a viscous friction contribution. Notice that the tensions in the cables must remain positive ($\Gamma_k \geq 0$, $k = 1, 2, 3$). 


\subsection{Control}

\subsubsection{Control issue and approach}

The goal of the controller is two-fold: 
\begin{itemize}
 \item have the position $q$ follow a reference $q_r$,
 \item maintain a positive tension in the cables.
\end{itemize}
The considered approaches are based on a linear model obtained by differentiation of the nonlinear model at the center of the workspace
\begin{equation}
	M_0 \, \ddot{q} + C_0 \, \dot{q} + K_0 \, q = J_0^\T \, \Gamma.
	\label{linéarisation}
\end{equation}
It is assumed that position and speed of the effector are fully available without restriction on the sensor bandwidth. Naturally, when the position varies from this central location, the performance of the considered control laws degrades and instability may be reached. This is precisely what we are going to investigate in section IV. 

\subsubsection{PD controller}


The first control law is designed such that the closed-loop system fits a generalized 2nd order ODE
\begin{equation}
	\ddot{q} + a_1 \, \dot{q} + a_0 \, (q - q_r)=0,
	\label{réponse visée}
\end{equation}
where $a_1 = 2 \xi \omega_0$ and $a_0 = {\omega_0}^2$, with $\omega_0$ the desired natural angular frequency (in rad/s$^{-1}$) and $\xi$ the desired damping.
 
The positioning issue is obtained by eliminating $\ddot{q}$ in (\ref{linéarisation}) and (\ref{réponse visée}), resulting in
\begin{equation}
	(C_0 - a_1 M_0) \, \dot{q} + (K_0 - a_0 M_0) \, q + a_0 M_0 \, q_r = J_0^\T \, \Gamma.
	\label{control}
\end{equation}
The solution can be written $\Gamma = \Gamma_p + \Gamma_t$ where $\Gamma_p = J_0^{\T \dag} \, \left( (C_0 - a_1 M_0) \, \dot{q} + (K - a_0 M_0) \, q \right)$ (the $\dag$ symbol stands for the Moore-Penrose pseudo-inverse) is the solution of the system of minimal norm and $\Gamma_t$ is chosen in the kernel of $J_0^\T$ so that it does not affect the movement. The solution $\Gamma_t = C_{moy} \begin{bmatrix}1&1&1 \end{bmatrix}^\T$ allows to set the average tension among the three cables at $C_{moy}$. For more details on tension management in CDPR, see for instance \cite{OA05}. 


The resulting control law includes a proportional-derivative feedback $\Gamma = K_d \, \dot{q} + K_p \, q + K_r \, q_r + C_{moy} \begin{bmatrix}1&1&1 \end{bmatrix}^\T$, where $K_d = J_0^{\T \dag} \, (C_0 - M_0 \bar{C})$, $K_p = J_0^{\T \dag} \, (K_0 - M_0 \bar{K})$ and $K_r = J_0^{\T \dag} \, M_0 \bar{K}$.

\begin{figure*}
\begin{center}
\begin{footnotesize}
 \begin{tabular}{ccc}
  \includegraphics[width=5.5cm]{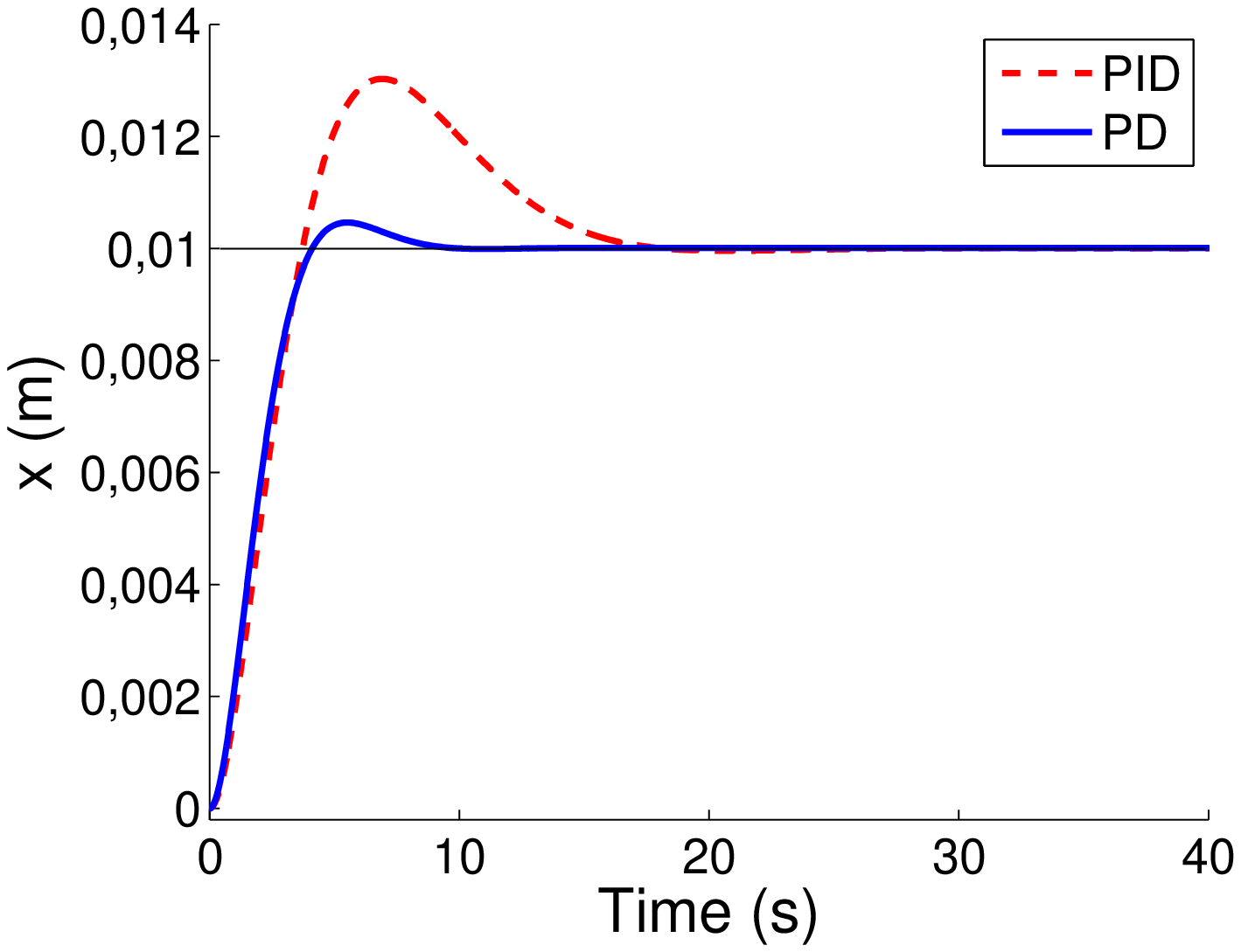} &
  \includegraphics[width=5.5cm]{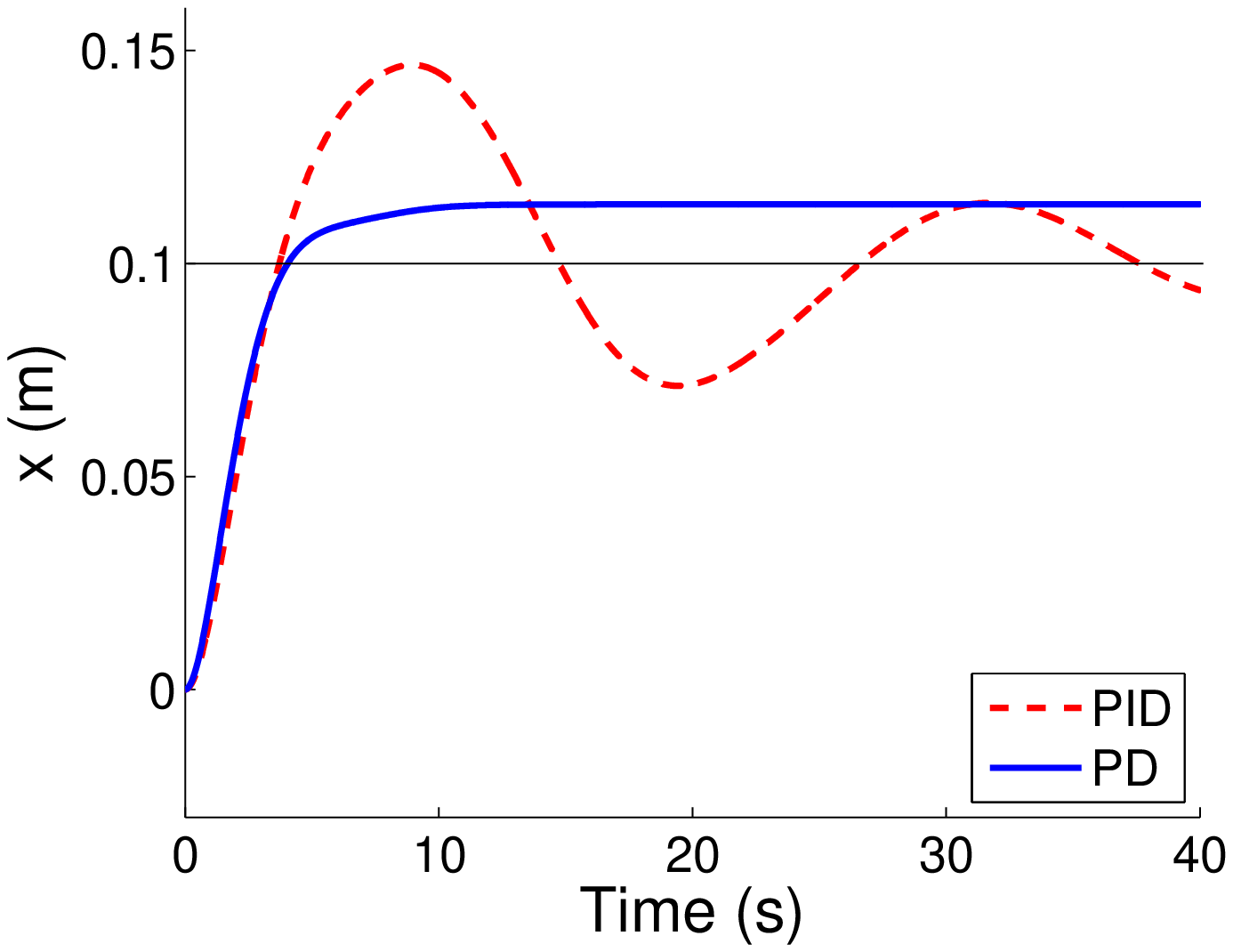} &
  \includegraphics[width=5.5cm]{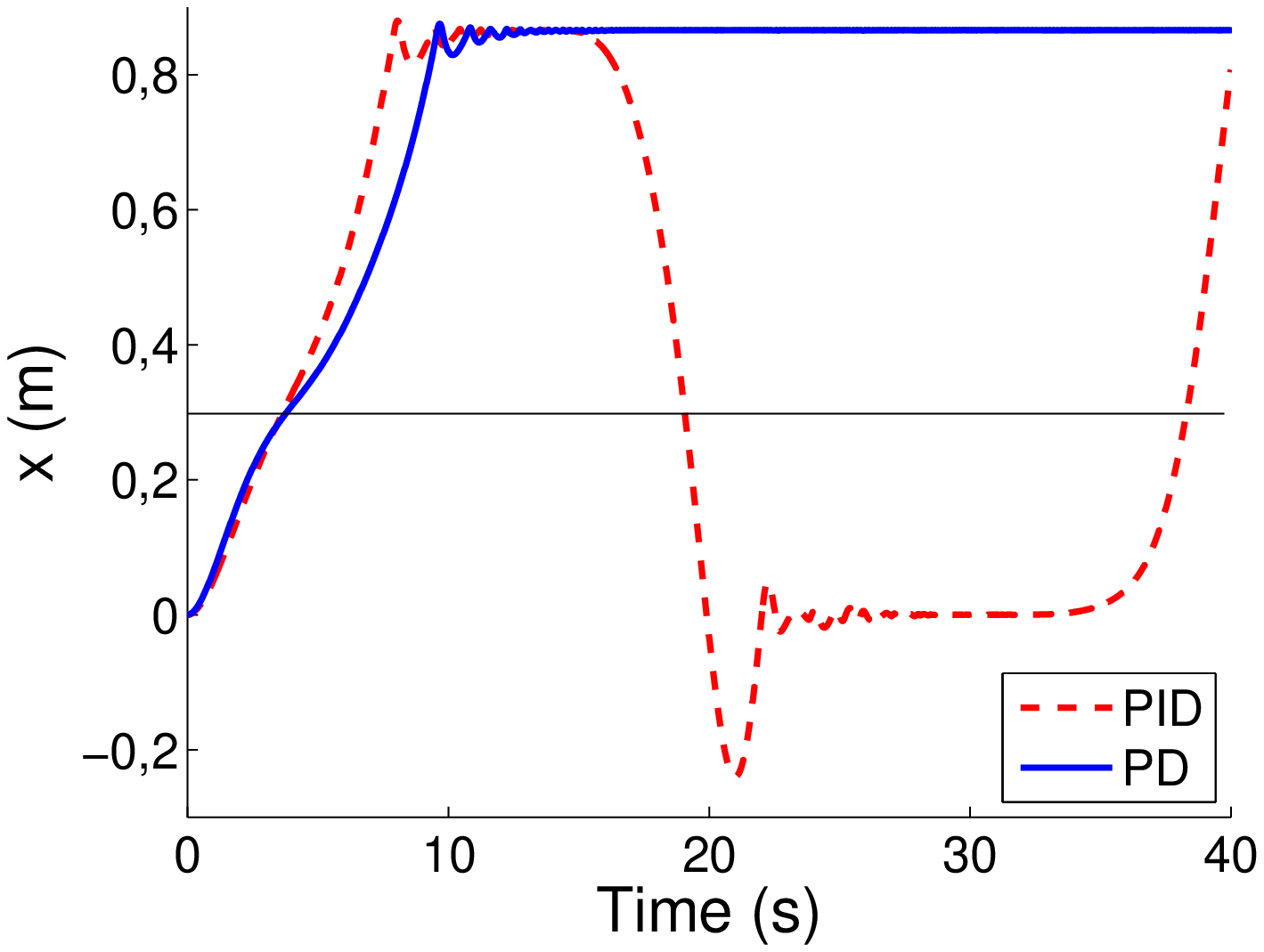} \\
  (a) $x_r$ = 0.01 m &
  (b) $x_r$ = 0.1 m &
  (c) $x_r$ = 0.3 m
 \end{tabular}
\end{footnotesize}
    \caption{Time response to a step of different amplitudes $x_r$ in the $x$ direction  (time-domain computations)} \label{fig:TimeResp}
\end{center}
\end{figure*}

\subsubsection{PID controller}

To ensure that the position reaches the reference in presence of modeling errors or constant perturbations, it is sensible to include some integral effect in the control law, thus introducing $I$ defined by $\dot{I} = q - q_r$ and considering $\Gamma_p = K_d \, \dot{q} + K_p\,  q + K_r \, q_r + K_i \, I$. 

The controller can be tuned so that the system behavior matches a reference system
\begin{equation}
	\dddot{q} + a_2 \, \ddot{q} + a_1 \, \dot{q} + a_0 \, (q - q_r) = 0,
	\label{réponse visée 2}
\end{equation}
with $K_d = J_0^{\T \dag} (C_0 - a_1 M_0)$, $K_p = J_0^{\T \dag} (K_0 - a_1 M_0)$ and $K_i = J_0^{\T \dag} a_0 J_0^{\T \dag} M_0$. 

The following pole placement strategy has been chosen: two complex poles corresponding to natural frequency of $\omega_0$ and a damping of $\xi$, plus a real pole $-\gamma$. For simplicity, the choice $\gamma = \xi \omega_0$ is done, i.e. the three poles have the same real part. The coefficients in (\ref{réponse visée 2}) are computed with $a_2 = \gamma + 2 \xi \omega_0$, $a_1 = \omega_0 \,(2 \xi \gamma + \omega_0^2)$ and $a_0 = \gamma \omega_0^2$.

\subsection{Simulation results}


The winders form an equilateral triangle with vertices located at $(\sqrt{3}/2,1/2)$, $(-\sqrt{3}/2,1/2)$ and $(0,-1)$.  The parameters of the model and the control laws are reported in Table~\ref{conf init}.
\begin{table}
	\caption{Dynamical parameters of the controlled systems}
	\label{conf init}
	\centering
	\begin{tabular}{|c||c|}
	\hline 
	 Parameter & Value \\
	\hline
	 $m$ & 1 kg        \\\hline
	 $I$ & 50.10$^{-6}$ kg.m$^2$ \\\hline
	 $R$ & 0.05 m         \\\hline
         $b$ & 0.7 N/(m/s)    \\\hline
         $\xi$ & 0.7        \\\hline
         $C_{moy}$ & 0.05 N.m \\
	\hline
	\end{tabular}
\end{table}

Time-domain simulations have been driven in order to evidence the behavior of the system with the proposed control laws. These consist in different steps of various amplitudes on reference $x_r$. The results for the different amplitudes with the PD and PID controllers are given in Fig.~\ref{fig:TimeResp}. For small amplitude steps (Fig.~\ref{fig:TimeResp}.a), both controllers provide a satisfying response and the system reaches the reference position with good accuracy. For medium amplitude (Fig.~\ref{fig:TimeResp}.b), the responses differ between the two controllers. One can notice that the response with the PID controller is much more oscillating while the response with the PD controller exhibits a static error. For large amplitude reference (Fig.~\ref{fig:TimeResp}.c), both controllers can be considered as unstable: the platform is attracted by the winder position. With the PID controller, the continuously integrating term results in oscillations between two winder positions. 



\section{RESULTS}\label{sec:results}

Two different numerical experiments are proposed to evaluate the abilities of MatCont in the context of this CDPR.  Bifurcation and stability analysis with respect to the reference positions is proposed in subsection \ref{sec:stab} for the PID and the PD controllers described in Section \ref{sec:CDPR}. Stability domains with respect to the variation of the natural frequency $\omega_0$ are plotted in subsection \ref{sec:stabilitydomain}. The continuation parameters used in MatCont are reported in Table \ref{tab:cont parameters}.
\begin{table}
	\caption{MatCont continuation parameters}
	\label{tab:cont parameters}
	\centering
	\begin{tabular}{|c||c||c|}
		\hline
		Parameter & Equilibrium  & Hopf curve \\
		          & (First step) & (Limit curve) \\
		\hline
		Number of points & 800 & 300 \\\hline
		Init Step Size & 0.001 & 0.01 \\\hline
		Max Step Size & 0.001 & 0.1\\\hline
		Min Step Size & $10^{-12}$ & $10^{-12}$ \\\hline
		Newton Tolerance & $10^{-6}$ & $10^{-6}$\\\hline
		Test Tolerance & $10^{-5}$ & $10^{-5}$ \\
		\hline
	\end{tabular}
\end{table}

\subsection{Bifurcation and stability analysis}\label{sec:stab}
As a first experiment, one-codim continuations are carried out from the stable equilibrium $(x,y,\dot{x},\dot{y})=(0,0,0,0)$ of the two CDPR models by varying $x_r$ along the horizontal axis, then by varying $y_r$ along the vertical symmetry axis as far as possible (see Fig.~\ref{fig:PIDstab}.a for the PID and Fig.~\ref{fig:PDstab}.a for the PD). This process enables to compute branches of equilibriums for the ODE system under study and to detect possible bifurcation points (limit points (LP), Hopf points (H) or branch points (BP)). As discussed in \ref{sec:MatCont}, these are indicators of stability loss. The solutions comprised between the stable equilibrium (0,0) and the first encountered  bifurcation point are stable ones, the others are unstable ones.  Two-codim continuations are then carried out from these bifurcation points to exhibit a fold curve (curve of limit points) assuming the PD controller, or a Hopf curve assuming the PID controller. Computations with the PID and the PD controllers are performed with $\omega_0=0.5$ rad.s$^{-1}$ and $\omega_0=0.8$ rad.s$^{-1}$, respectively, to get stable regions of similar area.

\subsubsection{PID controller}\label{sec:stabPID}
\begin{figure}
\begin{center}
\begin{footnotesize}
\begin{tabular}{cc}
\includegraphics[width=5cm]{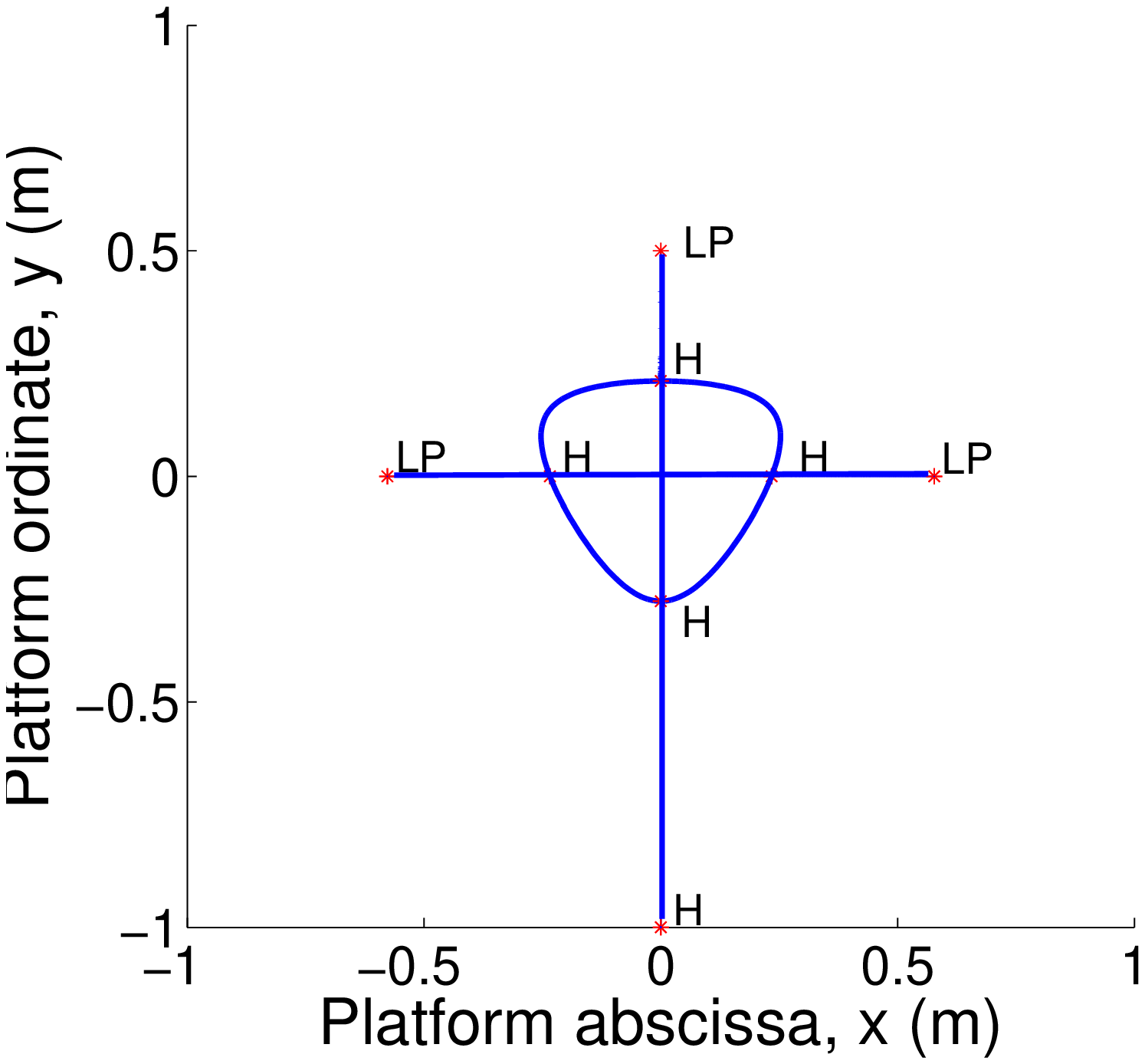} & 
\includegraphics[width=5.6cm]{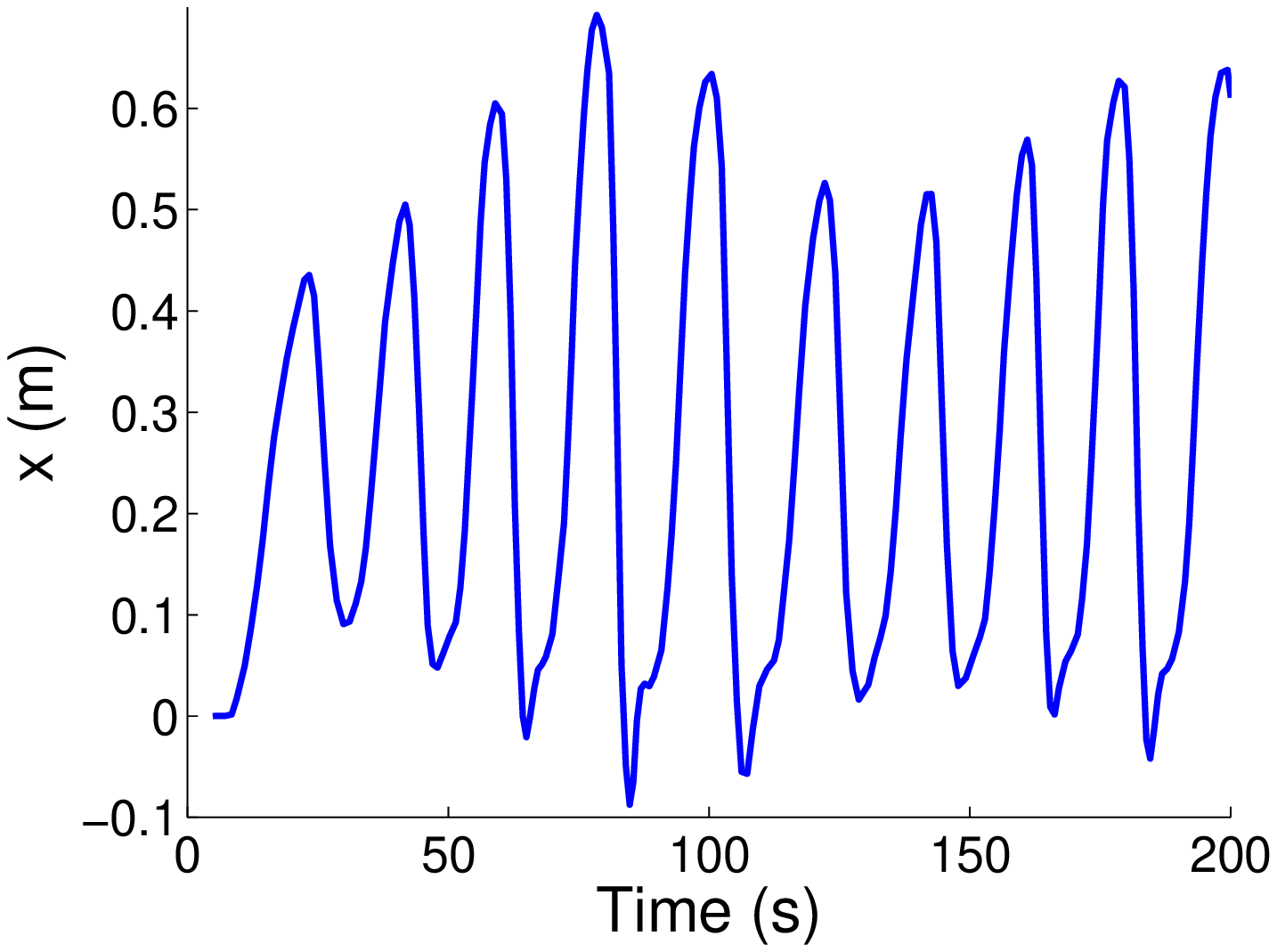} \\
(a) & (b)
\end{tabular}
    \caption{With the PID controller for $\omega_0=0.5$ rad.s$^{-1}$. (a) Bifurcation diagram. (b) Time-domain simulation at the Hopf point $(0.260,-0.010)$} \label{fig:PIDstab}
\end{footnotesize}
\end{center}
\end{figure}

The bifurcation diagram projected on the $(x,y)$ plane is plotted in Fig. \ref{fig:PIDstab}.a. As expected, the figure is symmetric with respect to the $y$-axis. MatCont detects a Hopf bifurcation at $q_r=(0.234,0)$ and a limit point at $q_r=(0.577,0)$ along the dynamical equilibrium curve computed by the one-codim continuation with respect to the positive coordinate $x_r$. The same observation holds for the other three semi-axes, except that a Hopf point is detected at the winder $A_3$. The limit points are located on the equilateral triangle defined by the winders. The branches linking (0,0) to the closest Hopf points are stable ones, while the branches linking these Hopf points to their corresponding limit/Hopf points are unstable ones. A 2-codim continuation in $x_r$ and $y_r$ is carried out from the Hopf point $(0.234,0)$. This results in a Hopf curve that links the already identified Hopf points and bounds the stability region. 
The Hopf point $(0.234,0)$ corresponds to an equilibrium point where two conjugate complex poles (eigenvalues of the Jacobian) become imaginary numbers. The robot reaches a stability limit and starts to oscillate. A simulation for this value (Fig \ref{fig:PIDstab}.b) confirms that the effector is submitted to periodic oscillations that are amplified in the unstable region.  

\subsubsection{PD controller}\label{sec:stabPD}
\begin{figure}
\begin{center}
\begin{footnotesize}
\begin{tabular}{cc}
\includegraphics[width=5cm]{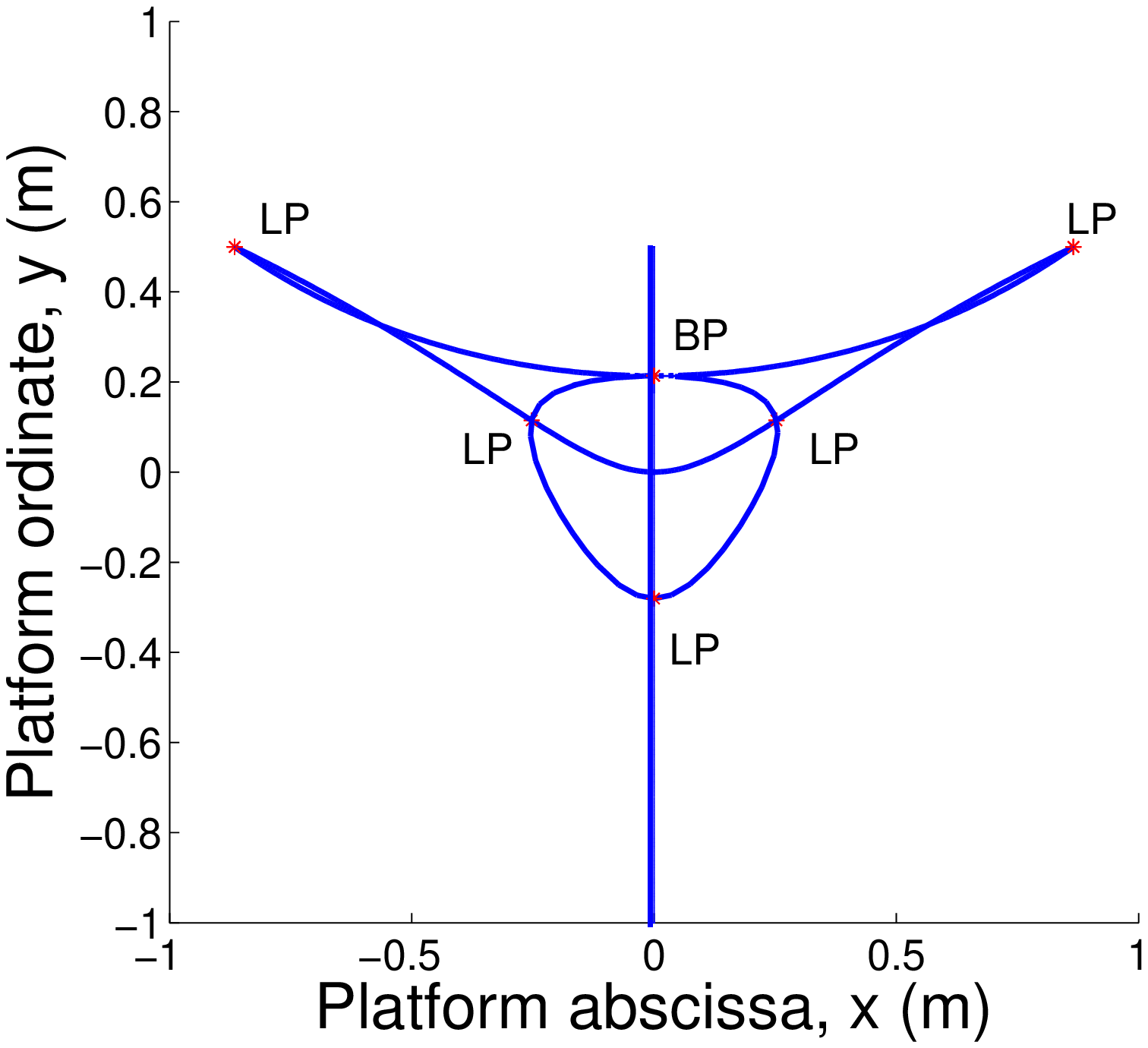} &
\includegraphics[width=6cm]{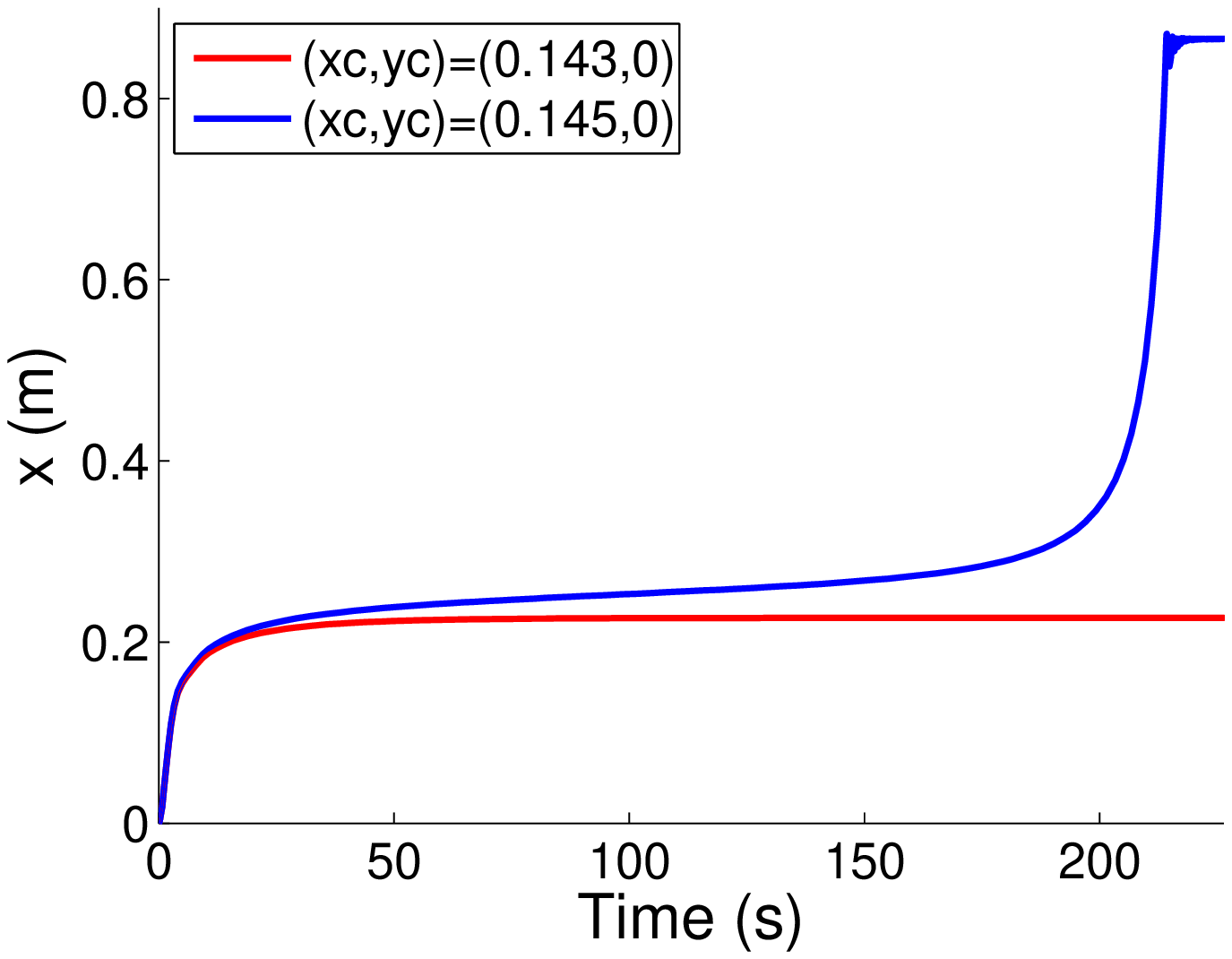}\\
(a) & (b)
\end{tabular}
\end{footnotesize}
\end{center}
    \caption{With the PD controller for $\omega_0=0.8$ rad.s$^{-1}$. (a) Bifurcation diagram. (b) Time-domain simulation close to the limit point $(0.144,0)$} \label{fig:PDstab}
\end{figure}

The bifurcation diagram projected on the $(x,y)$ plane is plotted on Fig.\ref{fig:PDstab}.a. Along the one-codim continuation with respect to positive $x_r$, MatCont detects a limit point for $(x_r,y_r)=(0.144,0)$. This corresponds to $(x,y)=(0.250,0.126)$ in the diagram. As expected, one observes that the computed branch is not a straight line in the $(x,y)$-plane due to the static error inherent to the PD controller. A branch point (BP) is detected along the one-codim continuation with respect to positive $y_r$. At this point, the equilibrium curve intersects with two unstable branches running to the two higher winders. A small perturbation applied close to the branch point drives the effector to one of the winders. The outer limit points correspond to winders $A_1$ and $A_2$. The two-codim continuation in $x_r$ and $y_r$ carried out from the limit point $(x_r,y_r)=(0.144,0)$  results in a fold curve that links the inner limit points and bounds the stability region.
Two time-domain simulations, Fig. \ref{fig:PDstab}.b, are performed close to the limit point $(x_r,y_r)=(0.144,0)$ to analyze the behavior of the mechanism. At this particular equilibrium, a pole is equal to $0$. As expected, the system is stable for (0.143, 0) and unstable for (0.145, 0) (the effector moves quickly towards the nearest winder), thus confirming the accuracy of the result provided by MatCont that the boundary of the stability domain lies between these two points. 

\subsection{Stability domains}\label{sec:stabilitydomain}
\begin{figure}
\begin{center}
    \includegraphics[width=8cm]{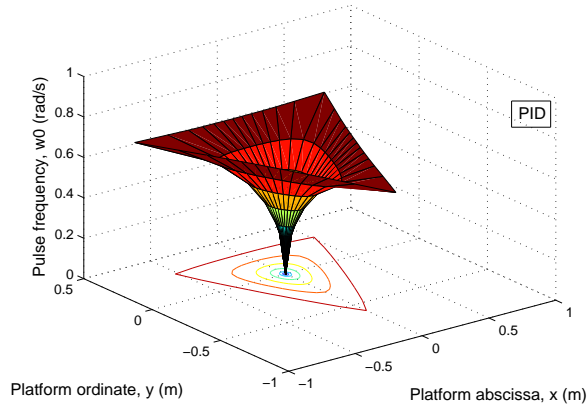}
\end{center}
    \caption{Stability domain for the PID controller}
		\label{fig:PIDisoline3D}
\end{figure}

\begin{figure}
\begin{center}
    \includegraphics[width=8cm]{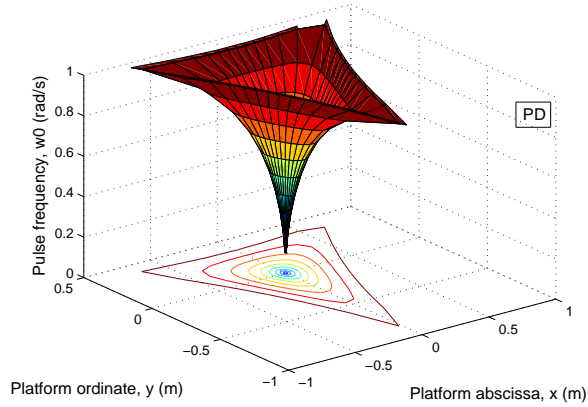}
    \caption{Stability domain for the PD controller}
		\label{fig:PDisoline3D}
\end{center}
\end{figure}

The desired natural frequency $\omega_0$ is now to be varied. The continuation process described in paragraph \ref{sec:stabPID} is carried out for different values of $\omega_0$. Experiments have been run using command lines together with CL\_MatCont. The stability curves (fold or Hopf) computed with respect to $(r_r,\phi_r)$ for given $\omega_0$ or with respect to $(r_r,\omega_0)$ for given angles $\phi_r$ may be used to build a 3D representation of the stability region, Fig. \ref{fig:PIDisoline3D} and \ref{fig:PDisoline3D}. Notice that MatCont is able to follow the curves up to  the lower value of $\omega_0$ using a two-codim continuation according to ($\omega_0$, $r_r$). At point ($x$,$y$) = (0,0), it detects a double Hopf point (HH) for the PID controller and a cusp point (CP) for the PD controller. This enables to calculate accurately the lower zone of the stability domains. 

The stability curves are assembled to figure out the stability domain boundaries as isolines of $\omega_0$ plotted in the $(x,y)$-plane. Boundary surfaces are deduced from a change of variables from Cartesian coordinates to cylindrical coordinates ($x_r = r_r \cos(\phi_r)$ and $x_r = r_r \sin(\phi_r)$) and an interpolation. The objective is to form a rectangular matrix and plot it using the Matlab's {\sl surfc} function, resulting in the stability domains presented in Figs. \ref{fig:PIDisoline3D} and \ref{fig:PDisoline3D}. 

\paragraph{PID controller}\label{sec:Workspace PID}
Hopf curves were computed for values of $\omega_0$ ranging from 0.2 to 0.65. The resulting stability domain is plotted in Fig. \ref{fig:PIDisoline3D}. 
The Hopf curves bound the inner stability region. The larger $\omega_0$, the larger the stable region. One explanation is that the real part of the poles at the nominal position, $-2 \xi \omega_0$, decreases with $\omega_0$ as $\xi$ is constant. For higher values of $\omega_0$, the stability domain spans the triangular workspace and is independent from $\omega_0$.  MatCont allows to compute the reference workspace $\mathcal{W}(x_r,y_r)$ of admissible $(x_r,y_r)$ as well. This is identical to  $\mathcal{W}(x,y)$ as expected from the use of a PID controller. 

The study did not determine upper limits as no bandwidth restrictions have been introduced in the actuators nor in the sensors; and no delay is considered in the control. 

\begin{figure}
\begin{center}
\begin{footnotesize}
\hfil\includegraphics[width=5cm]{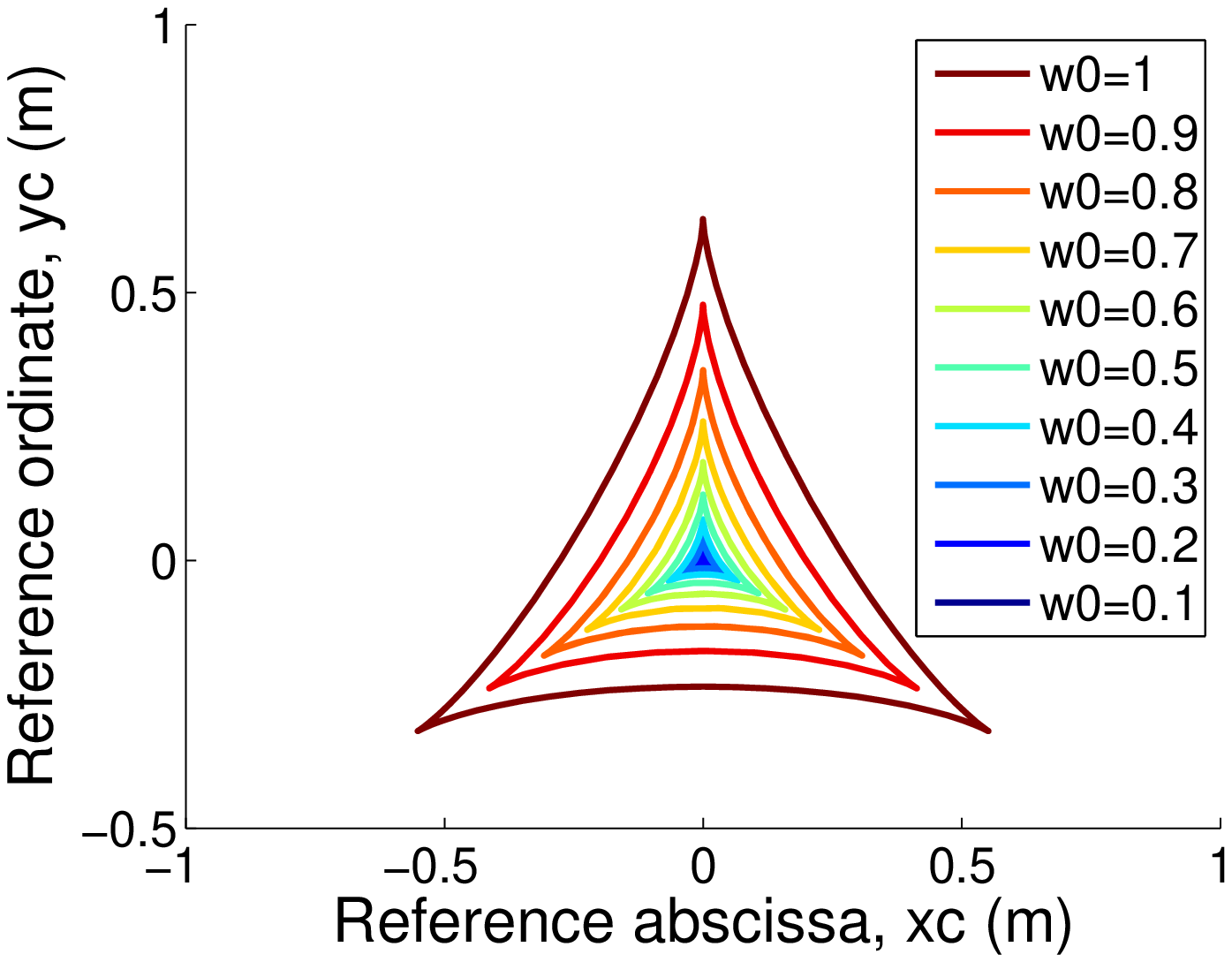}\hfil~\hfil \includegraphics[width=5cm]{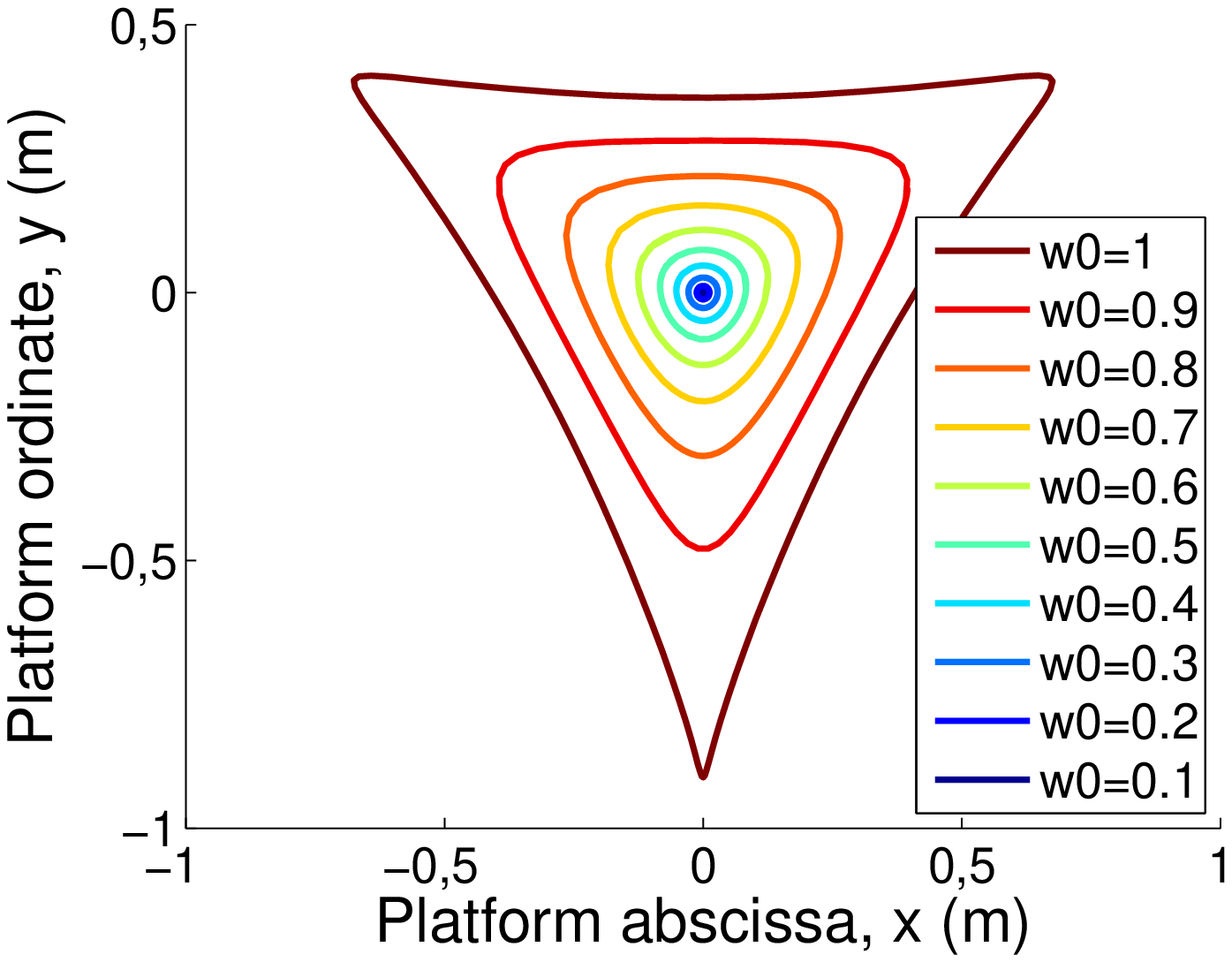}\hfil\\
\hfil(a) ($x_r,y_r$) plane \hfil~\hfil(b) ($x$, $y$) plane\hfil
\end{footnotesize}
\end{center}
    \caption{Stability domains for the PD controller. (a): position reference, (b): platform position} \label{fig:PDstab2}
\end{figure}

\paragraph{PD controller}
Figs. \ref{fig:PDisoline3D} and \ref{fig:PDstab2} show the fold curves computed for values of $\omega_0$ ranging from 0.1 and 1. together with the reconstructed surface. The computed curves show again that the workspace area grows when $\omega_0$ increases. By comparing Fig.~\ref{fig:PDisoline3D}.a and Fig.~\ref{fig:PDisoline3D}.b, one can see that the stability domains in the ($x_r$, $y_r$) domain or in the ($x$, $y$) plane do not have the same shape. As observed in simulation, without integral term, the platform is likely to move towards the nearest winder rather than stay at its reference position. Notice that the provided curves correspond to the limit of the stability domain. Inside the stability domain, this surprising tendency is weaker and satisfying tracking performance can be obtained in the neighborhood of the nominal position. 

\section{CONCLUSIONS}

In this paper, the use of MatCont -- a numerical package dedicated to continuation analysis -- has been investigated for the determination of the stability domain. The considered case-study is a simplistic cable-driven parallel robot with only two degrees of freedom. However, its nonlinear model exhibits non-rational nonlinearities, which make it a not-so-easy target for usual robustness analysis tools. With a very limited number of manipulations and a reasonable computational time, MatCont determines limits of the stability domain and allows to show the influence of the controller tuning. These positive results make MatCont a potential tool of control engineers and researchers.



\end{document}

%% file: robot_pspdf.tex
\begin{picture}(0,0)%
\includegraphics{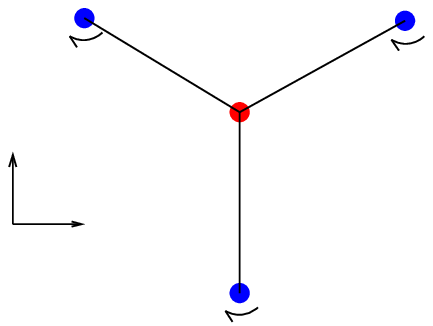}%
\end{picture}%
\setlength{\unitlength}{4144sp}%
\begingroup\makeatletter\ifx\SetFigFont\undefined%
\gdef\SetFigFont#1#2#3#4#5{%
  \reset@font\fontsize{#1}{#2pt}%
  \fontfamily{#3}\fontseries{#4}\fontshape{#5}%
  \selectfont}%
\fi\endgroup%
\begin{picture}(2135,1630)(1111,-3532)
\put(1440,-3214){\makebox(0,0)[lb]{\smash{{\SetFigFont{8}{9.6}{\rmdefault}{\mddefault}{\updefault}{\color[rgb]{0,0,0}$x$}%
}}}}
\put(1126,-2966){\makebox(0,0)[lb]{\smash{{\SetFigFont{8}{9.6}{\rmdefault}{\mddefault}{\updefault}{\color[rgb]{0,0,0}$y$}%
}}}}
\put(2461,-2670){\makebox(0,0)[lb]{\smash{{\SetFigFont{7}{8.4}{\rmdefault}{\mddefault}{\updefault}{\color[rgb]{0,0,0}$E(x,y)$}%
}}}}
\put(2444,-3315){\makebox(0,0)[lb]{\smash{{\SetFigFont{8}{9.6}{\rmdefault}{\mddefault}{\updefault}{\color[rgb]{0,0,0}$A_3(x_3,y_3)$}%
}}}}
\put(3122,-2372){\makebox(0,0)[lb]{\smash{{\SetFigFont{8}{9.6}{\rmdefault}{\mddefault}{\updefault}{\color[rgb]{0,0,0}$\theta_1$}%
}}}}
\put(2907,-2058){\makebox(0,0)[lb]{\smash{{\SetFigFont{8}{9.6}{\rmdefault}{\mddefault}{\updefault}{\color[rgb]{0,0,0}$A_1(x_1,y_1)$}%
}}}}
\put(1650,-2355){\makebox(0,0)[lb]{\smash{{\SetFigFont{8}{9.6}{\rmdefault}{\mddefault}{\updefault}{\color[rgb]{0,0,0}$\theta_2$}%
}}}}
\put(1452,-2025){\makebox(0,0)[lb]{\smash{{\SetFigFont{8}{9.6}{\rmdefault}{\mddefault}{\updefault}{\color[rgb]{0,0,0}$A_2(x_2,y_2)$}%
}}}}
\put(2206,-3481){\makebox(0,0)[lb]{\smash{{\SetFigFont{8}{9.6}{\rmdefault}{\mddefault}{\updefault}{\color[rgb]{0,0,0}$\theta_3$}%
}}}}
\end{picture}%

%% file: StabCont.bbl
\begin{thebibliography}{99}
\bibitem{AG90} E. L. Allgower, and K. Georg, {\it Numerical Continuation Methods: an Introduction}, New-York: Springer-Verlag, 1990, pp. 390.  
\bibitem{BP13}
T. Bruckmann and A. Pott, {\it Proceedings of the First International Conference on Cable-Driven Parallel Robots}, Springer, 2013
\bibitem{BP15}
T. Bruckmann and A. Pott, {\it Proceedings of the Second International Conference on Cable-Driven Parallel Robots}, Springer, 2015
\bibitem{charpentier12} I. Charpentier, On higher-order differentiation in nonlinear mechanics, {\it Optim. Method. Softw.}, vol. 27, no. 2, pp.  221--232, Feb. 2012.
\bibitem{CDPF07} B. Cochelin, N. Damil, and M. Potier-Ferry, M., {\it M\'{e}thode asymptotique num\'{e}rique}, Paris: Herm\`es-Lavoisier, 2007, pp.298.
\bibitem{MatCont} A. Dhooge, W. Govaerts, Yu. A. Kuznetsov, W. Mestrom, A. M. Riet, and B. Sautois, 
        {MatCont and CL MatCont: Continuation toolboxes in matlab (MatCont manual)}, 2006. 
\bibitem{FAG96} E. Feron, P. Apkarian, and P. Gahinet, Analysis and synthesis of robust control systems via parameter-dependent Lyapunov functions, {\it IEEE Trans. Autom. Control.}, vol. 41, no. 7, pp. 1041--1046, Jul. 1996.
\bibitem{DKK91}E. Doedel, H. Keller, and J.-P. Kernevez, Numerical analysis and control of bifurcation problems (i) bifurcation in finite dimensions, {\it Int. J. Bifurcat. Chaos}, vol. 1, no. 4, pp. 493--520, Sep. 1991.
\bibitem{haug96} E. Haug, C. Luh, F. Adkins, and J. Wang, Numerical algorithms for mapping boundaries of manipulator workspaces, {\it J. Mech. Des.-T. ASME}, vol. 118, no. 2, pp. 228--234, Jun. 1996.
\bibitem{crmeca}G. Hentz, I. Charpentier, and P. Renaud, Higher-order continuation for the determination of robot workspace boundaries, {\it CR Meca.}, vol. 344, p. 95--101, 2016.
\bibitem{Kel77} H. B. Keller, Numerical solution of bifurcation and nonlinear eigenvalue problems, in {\it Applications of Bifurcation Theory}, P. Rabinowitz, Ed., New-York: Academic Press, 1977, pp. 359--384.
\bibitem{Kel87} H. B. Keller, {\it Lectures on numerical methods in bifurcation problems}, Berlin: Springer-Verlag, 1987, pp. 155.
\bibitem{KZW06} K. Kozak, Q. Zhou, and J. Wang, Static analysis of cable-driven manipulators with non-negligible cable mass, {\em IEEE Trans. Robotics}, 
vol. 22, no. 3, pp. 425 - 433, 2006.
\bibitem{LCCG12} E. Laroche, R. Chellal, L. Cuvillon, and J. Gangloff, A preliminary study for H$_\infty$ control of parallel cable-driven manipulators, {\it First International Conference on Cable-Driven Parallel Robots}, Stuttgart, Aug. 2012.
\bibitem{litvin80} F. L. Litvin, Application of theorem of implicit function system existence for analysis and synthesis of linkages, {\it Mech. Mach. Theory}, vol. 15, no. 2, pp. 115--125, 1980.
\bibitem{OA05} S. R. Oh and S. Agrawal, Cable suspended planar robots with redundant cables: controllers with positive tensions, {\it IEEE Trans. Robot.}, vol. 21, no. 3, pp. 457-465, June 2005.
\bibitem{Packard93} A. Packard, and J. Doyle, The complex structured singular value, {\it Automatica}, vol. 29, no. 1, pp. 71--109, Jan. 1993.
\bibitem{Sey79} R. Seydel, Numerical computation of branch points in nonlinear equations, {\it Numer. Math.}, vol. 33, no. 3, pp. 339--352, Sep. 1979.
\bibitem{Sey09} R. Seydel, {\it Practical Bifurcation and Stability Analysis}, 3rd edition, New York: Springer, 2009, pp. 483. 

\end{thebibliography}
